\documentclass[11pt]{article}
\usepackage{amssymb,amsmath,enumerate,latexsym, graphicx, epsfig}
\usepackage[latin1]{inputenc}

\newtheorem{theorem}{Theorem}

\newtheorem{conjecture}{Conjecture}
\newtheorem{definition}{Definition}

\newenvironment{proof}{\smallskip\noindent{\it Proof.}\hskip \labelsep}
                        {\hfill\penalty10000\raisebox{-.09em}{$\Box$}\par\medskip}

\def\R{\mathbb{R}}

\def\N {{\Bbb N}}
\def\T {{\Bbb T}}
\def\S {{\Bbb S}}
\def\Z {{\Bbb Z}}

\begin{document}

\title{Minimal surfaces with genus zero}%
\author{M. Magdalena Rodr\'\i guez}%
\date{\mbox{}}
\maketitle

\noindent {\sc Abstract.} {\footnotesize A very interesting problem in
  the classical theory of minimal surfaces consists of the
  classification of such surfaces under some geometrical and
  topological constraints. In this short paper, we give a brief
  summary of the known classification results for properly embedded
  minimal surfaces with genus zero in $\R^3$ or quotients of $\R^3$ by
  one or two independent translations.  This does not intend to be an
  exhaustive review of the tools or proofs in the field, but a simple
  explanation of the currently known results.  }

\section{Introduction}

The classical theory of minimal surfaces, whose roots go back to
Euler and Lagrange in the 18th century, still remains extremely active
nowadays.  New examples discovered in an explosion of activity in the
eighties have gradually focused the subject on the problem of
classification.  Recently, important advances have been made
towards the goal of characterizing the minimal surfaces under some
topological and/or geometrical constraints, and new families of
unexpected examples have been found.

In this short paper, we briefly expose the history and some known
classification results for minimal surfaces whose genus is zero (i.e.
those which are topologically a punctured sphere).  For more detailed
reviews, see~\cite{cmCourant,cm34, mpe1} and the references therein.

\section{Genus zero minimal surfaces in $\R^3$}\label{sec1}
\begin{center}
  \begin{figure}
    \includegraphics[width=12.5cm,height=3.2cm]{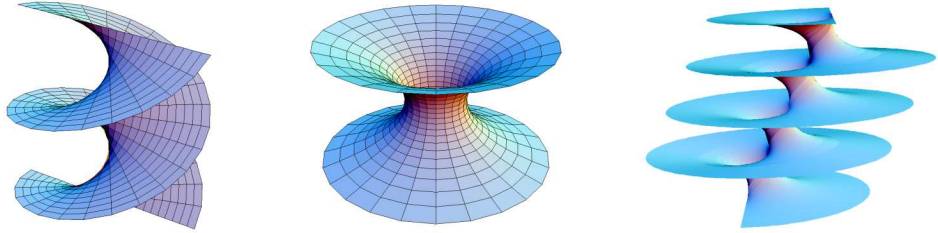}
    \caption{The helicoid (left), the catenoid (center) and a Riemann
      minimal example (right).}
    \label{sup1}
  \end{figure}
\end{center}

In 1740, Euler constructed the first nonplanar minimal surface: the
{\it catenoid} (see Figure~\ref{sup1}, center).  He also showed that
this is the only nonplanar, complete minimal surface of revolution in
$\R^3$. The following theorem provides us another important
characterization of the catenoid.
\begin{theorem}[L\' opez \& Ros, \cite{lor1}]\label{loros}
  The only nonplanar properly embedded minimal surface in $\R^3$ with
  finite total curvature\footnote{Let $M\subset\R^3$ be a minimal
    surface, and $K$ its Gaussian curvature. We define the {\it total
      curvature} of $M$ as $C(M)=\int_M K$.} and genus zero is the
  catenoid.
\end{theorem}

Thirty years later the discovery of the catenoid, Meusnier found
another minimal example: the {\it helicoid} (see Figure~\ref{sup1},
left), that was classified in 1842 by Catalan as the only nonplanar
ruled minimal surface in $\R^3$.  Recently, the helicoid has also been
classified, in a beautiful application of the outstanding theory
developed by Colding and Minicozzi, as the only properly embedded
minimal surface with genus zero and one end.
\begin{theorem}[Meeks \& Rosenberg, \cite{mr8}]\label{hel}
  A properly embedded simply - connected minimal surface in $\R^3$ is
  either the plane or the helicoid.
\end{theorem}

A theorem by Collin~\cite{col1} assures that, if $M\subset\R^3$ is a
properly embedded minimal surface with at least two ends, then $M$ has
finite topology if and only if it has finite total curvature. Hence,
we obtain from Theorems~\ref{loros} and~\ref{hel} that the only
properly embedded minimal surfaces in $\R^3$ with genus zero and
finite topology (i.e. with finitely many ends) are the plane, the
helicoid and the catenoid.  In particular, there are no properly
embedded minimal surfaces with genus zero and $n$ ends, for any
$n\in\N, n\geq 3$.  However, there exist properly embedded minimal
surfaces with genus zero and infinitely many ends: the {\it Riemann
  minimal examples}.  In a posthumously published paper,
Riemann 
constructed a 1-paramenter family ${\cal R}=\{R_t\}_{t>0}$ of minimal
surfaces in $\R^3$, and classified them together with the plane, the
helicoid and the catenoid as the only minimal surfaces of $\R^3$
foliated by straight lines and/or circles.  Each $R_t$ has genus zero
with infinitely many ends in $\R^3$, but its geometry simplifies
considerably modulo its symmetries: $R_t$ is invariant by the
translation by $(t,0,2)$ and has two planar ends and genus one in the
quotient by $(t,0,2)$. Viewed in $\R^3$, $R_t$ has two limit
ends\footnote{A limit end $e$ of a noncompact surface $M$ is an
  accumulation point of the set $E(M)$ of ends of $M$; this makes sense
  since $E(M)$ can be endowed with a natural topology which makes it a
  compact, totally disconnected subspace of the real interval $[0,1]$,
  see~\cite{mpe1}.}: one top and one bottom limit end (see
Figure~\ref{sup1}, right).  These examples have also been classified
as follows.

\begin{theorem}[Meeks, P\'erez \& Ros, \cite{mpr1}]
  The Riemann minimal examples are the unique periodic\footnote{A
    surface in $\R^3$ is said to be {\it singly, doubly or triply
      periodic} when it is invariant by a discrete infinite group $G$
    of isometries of $\R^3$ of rank one, two or three, respectively,
    that acts properly and discontinuously. After passing to a
    finitely sheeted covering, we can assume that the flat 3-manifold
    $\R^3/G$ is either $\S^1\times\R^2, \R^3/S_\theta, \T^2\times\R$
    or $\T^3$, where $S_\theta$ denotes a screw motion symmetry and
    $\T^k$ is a $k$-dimensional flat torus, $k=2,3$.} non
  simply-connected properly embedded minimal surfaces in $\R^3$ with
  genus zero.
\end{theorem}

The extension of the above characterization by eliminating the
hypothesis on the periodicity constitutes the following conjecture, by
the same authors.

\begin{conjecture}[Genus zero conjecture, \cite{mpr3}]\label{conj}
  If $M\subset\R^3$ is a properly embedded minimal surface with genus
  zero and infinitely many ends, then $M$ is a Riemann minimal example.
\end{conjecture}

Meeks, P\'erez and Ros have proven some strong partial results in this
direction; for example, they prove in~\cite{mpr4} that a surface $M$
in the hypothesis of Conjecture~\ref{conj} cannot have exactly one
limit end.

\section{Genus zero minimal surfaces in $\T^3$ }\label{sec4}
Let us denote by $\T^3$ a 3-dimensional flat torus obtained as a
quotient of $\R^3$ by three independent translations.  Since $\T^3$ is
a compact space, every minimal surface $M\subset\T^3$ is also compact,
so it cannot have any end.  In particular, a genus zero minimal
surface in $T^3$ is topologically a sphere, and so it lifts to a minimal
sphere in $\R^3$, which is impossible.
Hence there are no properly embedded, triply periodic, minimal
surfaces in $\R^3$ with genus zero in the quotient.

\section{Genus zero minimal surfaces in $\T^2\times\R$}\label{sec3}
In 1835, Scherk 
showed a properly embedded minimal surface $S$ in $\R^3$
invariant by two independent translations. The surface $S$, which is
obtained from a minimal graph over the unit square with boundary data
$\pm\infty$ disposed alternately, is invariant by $(2,0,0), (0,2,0)$.
This surface $S$ was generalized later on to a 1-parameter family of
properly embedded minimal surfaces $S_\theta$ invariant by
translations of vectors $(2,0,0),(2\cos\theta,2\sin\theta,0)$, for
$\theta\in(0,\frac{\pi }{2}]$ (in particular, $S= S_\frac{\pi }{2}$).
These surfaces are known as {\it doubly periodic Scherk minimal
  examples} and have genus zero and four ends asymptotic to flat
vertical annuli in the quotient by such translations. Annular ends of
this kind are usually called {\it Scherk-type ends}.

\begin{theorem}[Lazard-Holly \& Meeks, \cite{lhm}]\label{2pScherk}
  The only nonplanar properly embedded minimal surfaces in
  $\T^2\times\R$ with genus zero are the doubly periodic Scherk
  minimal examples.
\end{theorem}

\section{Genus zero minimal surfaces in $\S^1\times\R^2$}\label{sec2}
The first known properly embedded minimal surfaces of $\R^3$ invariant
by only one independent translation were the conjugate surfaces
$S_\theta^*$ of the doubly periodic Scherk minimal examples,
$\theta\in(0,\frac{\pi }{2}]$, called {\it singly periodic Scherk
  minimal examples}. They may be viewed as the desingularization of
two vertical planes meeting at angle $\theta$, and are invariant by
the translation by $T=(0,0,2)$.  In the quotient by $T$, denoted as
$\S^1\times\R^2$, each $S_\theta$ has genus zero and four Scherk-type
ends.

H. Karcher~\cite{ka4} generalized the previous Scherk examples by
constructing, for each natural $n\geq 2$, a $(2n-3)$-parameter family
of properly embedded minimal surfaces in $\R^3$, invariant by $T$, and
with genus zero and $2n$ Scherk-type ends in the quotient. These
surfaces are called {\it saddle towers}.
Let us now recall their construction: Consider any convex polygonal
domain $\Omega_n$ whose boundary consists of $2n$ edges of length one,
with $n\geq 2$, and mark its edges alternately by $\pm\infty$.  Assume
$\Omega_n$ is non-special (see definition~\ref{def1} below).  By a
theorem of Jenkins and Serrin~\cite{jes1}, there exists a minimal
graph defined on $\Omega_n$ which diverges to $\pm\infty$, as
indicated by the marking, when we approach to the edges of $\Omega_n$.
The boundary of this minimal graph consists of $2n$ vertical lines
above the vertices of $\Omega_n$. Hence the conjugate minimal surface
of this graph is bounded by $2n$ horizontal symmetry curves, lying between
two horizontal planes separated by distance one.  By reflecting
about one of the two symmetry planes, we obtain a fundamental domain
for a saddle tower with $2n$ Scherk-type ends. 

\begin{definition}
\label{def1}
We say that a convex polygonal domain with $2n$ unitary edges is
special if $n\geq 3$ and its boundary is a parallelogram with two
sides of length one and two sides of length $n-1$.
\end{definition}

\begin{theorem}[P\'erez \& Traizet, \cite{PeTra1}]
  If $M\subset\S^1\times\R^2$ is a nonplanar, properly embedded
  minimal surface with genus zero and finitely many Scherk-type ends,
  then $M$ is a saddle tower.
\end{theorem}

Now we look for properly embedded minimal surfaces in $\S^1\times\R^2$
with genus zero and infinitely many ends. In a first step, and
following the arguments in Section~\ref{sec1}, we try to classify
those which are periodic. There are two families in this setting: the
singly periodic liftings of the 1-parameter family of doubly periodic
Scherk minimal examples and singly periodic liftings of the
3-parametric family of {\it KMR examples}~\cite{ka4,mrod1} (also
called toroidal halfplane layers).  These later examples have been
classified in~\cite{PeRoTra1} as the only properly embedded, doubly
periodic, minimal surfaces in $\R^3$ with parallel ends and genus one
in the quotient. If we consider the quotient of a KMR example only by
the period at its ends, we obtain a periodic minimal surface in
$\S^1\times\R^2$ with genus zero, infinitely many ends, and two limit
ends: one top and one bottom. 

\begin{theorem}\label{th}
  Let $M\subset\S^1\times\R^2$ be a (periodic) properly embedded minimal
  surface with genus zero so that its lifting to $\R^3$ is a doubly
  periodic surface. Then $M$ is either a doubly periodic Scherk
  minimal example (with one limit end) or a KMR example (with two
  limit ends).
\end{theorem}
\begin{proof}
  Let $M\subset\S^1\times\R^2$ be in the hypothesis of
  Theorem~\ref{th}. Then, its lifting $\widetilde M$ to $\R^3$ is a
  properly embedded minimal surface invariant by two independent
  translations $T_1,T_2$, which can only have genus zero or one in the
  quotient $\overline{M}$ by both translations.

  Suppose that $\widetilde M$ has non-parallel ends. Thus $T_1$ and
  $T_2$ are generated by the period at the ends of $\widetilde M$.
  When $\overline{M}$ has genus zero, it must be a doubly periodic
  Scherk minimal surface, by Theorem~\ref{2pScherk}.  Suppose
  $\overline{M}$ has genus one, and let $\{\gamma_1,\gamma_2\}$ be a
  homology basis of $H_1(M,\Z)$. By possibly adding to $\gamma_i$,
  $i=1,2$, a finite number of loops around the ends of $\overline M$,
  we obtain a homology basis
  $\{\widetilde\gamma_1,\widetilde\gamma_2\}$ of $H_1(M,\Z)$ so that
  the period of $\overline M$ along $\widetilde\gamma_i$ vanishes,
  which is not possible. Hence, if $\widetilde M$ has non-parallel
  ends, it must be a doubly periodic Scherk minimal surface.

  Finally, suppose that $\widetilde M$ has parallel ends. By
  Theorem~\ref{2pScherk}, $\overline{M}$ must have genus one, and so
  it is a KMR example,~\cite{PeRoTra1}. This finishes
  Theorem~\ref{th}.
\end{proof}

It was expected that there were no more properly embedded minimal
surfaces in $\S^1\times\R^2$ with genus zero and infinitely many ends
other than the periodic ones, as in the case of $\R^3$. But this is
not true, as the following theorem shows.

\begin{theorem}[Mazet, \_\_ \& Traizet, \cite{mrt}]\label{mrt}
  For each non special unbounded convex polygonal domain $\Omega$ with
  unitary edges (see Definition~\ref{def2}), there exists a non
  periodic, properly embedded minimal surface
  $M_\Omega\subset\S^1\times\R^2$ with genus zero, infinitely many
  ends and one limit end, whose conjugate surface can be obtained from
  a minimal graph over $\Omega$ with boundary values $\pm\infty$
  disposed alternately.
\end{theorem}

\begin{definition}
  \label{def2}
  An unbounded convex polygonal domain is said to be special when its
  boundary is made of two parallel half lines and one edge of length
  one (such a domain may be seen as a limit of special domains with
  $2n$ edges, when $n\to\infty$).
\end{definition}

\begin{figure}
  \begin{center}
    \includegraphics[width=10.93cm,height=5cm]{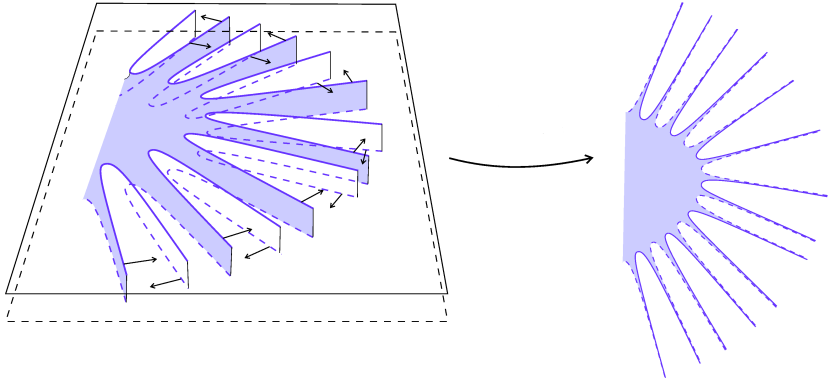}
    \caption{The sketch of half a surface $M_\Omega$ satisfying Theorem~\ref{mrt}.}
    \label{InfSaddle}
  \end{center}
\end{figure}

The surfaces $M_\Omega$ appearing in Theorem~\ref{mrt} are obtained by
taking limits of saddle towers $M_n$, each $M_n$ with $4n$ ends.
In~\cite{mrod3} we have proven that the other possible limits for a such
sequence $\{M_n\}_n$ of saddle towers are: the singly periodic Scherk
minimal example $S_{\frac{\pi}{2}}^*$, every doubly periodic Scherk
  minimal example, a KMR example $M_{\theta,\alpha,0}$ studied
  in~\cite{mrod1} or, after blowing up, a catenoid.

\section*{Acknowledgments}
The author would like to thank Joaqu\'\i n P\'erez for helpful conversations.

Research partially supported by grants from R\'egion
  Ile-de-France and a MEC/FEDER grant no.  MTM2004-02746.

\end{document}